\documentclass[12pt,reqno]{amsart}
\usepackage{amsmath, amsfonts, amssymb, amsthm, hyperref}
\usepackage{bm}
\allowdisplaybreaks[4]
\def\bg{\bigg}
\def\({\bg(}
\def\){\bg)}
\def\t{\text}
\def\f{\frac}

\def\sgn{{\rm sgn}}

\def\ls{\leqslant}

\def\Proof{\noindent{\it Proof}}

\def\Z{\mathbb Z}

\def\1{{\bf 1}}\def\c{{\bf c}}

\def\pmod #1{\ ({\rm{mod}}\ #1)}
\def\mod #1{\ {\rm mod}\ #1}

\def\<{\langle}
\def\>{\rangle}
\def\Ack{\medskip\noindent {\bf Acknowledgments}}
\theoremstyle{plain}
\newtheorem{theorem}{Theorem}[section]
\newtheorem{lemma}{Lemma}
\newtheorem{corollary}{Corollary}
\newtheorem{conjecture}{Conjecture}

\theoremstyle{definition}

\theoremstyle{remark}
\newtheorem{remark}{Remark}

\makeatletter
\@namedef{subjclassname@2020}{%
  \textup{2020} Mathematics Subject Classification}
\makeatother
 \vspace{4mm}

\begin{document}
\medskip

\title[Trinomial coefficients and related matrices over finite fields]
{Trinomial coefficients and related matrices over finite fields}
\author{Yue-Feng She}
\address {(Yue-Feng She) Department of Mathematics, Nanjing
	University, Nanjing 210093, People's Republic of China}
\email{she.math@smail.nju.edu.cn}

\author{Hai-Liang Wu}
\address {(Hai-Liang Wu) School of Science, Nanjing University of Posts and Telecommunications, Nanjing 210023, People's Republic of China}
\email{whl.math@smail.nju.edu.cn}

\keywords{Trinomial coefficients, Finite Fields, Determinants.
\newline \indent 2020 {\it Mathematics Subject Classification}. Primary 05A19, 11C20; Secondary 15A18, 15B57, 33B10.
\newline \indent This work was supported by the Natural Science Foundation of China (Grant no. 12101321).}

\begin{abstract}
In this paper, with the help of trinomial coefficients we study some arithmetic properties of certain determiants involving reciprocals of binary quadratic forms over finite fields.
\end{abstract}
\maketitle

\section{Introduction}
\setcounter{lemma}{0}
\setcounter{theorem}{0}
\setcounter{equation}{0}
\setcounter{conjecture}{0}
\setcounter{remark}{0}
\setcounter{corollary}{0}

\subsection{Background and Motivations} 
Determinants over finite fields have extensive applications in both number theory and combinatorics. Readers may refer to the survey papers \cite{K1,K2} for the recent progresses on this topic. Throughout this paper, for any prime $p$, the symbol $\mathbb{F}_p$ denotes the finite field with $p$ elements.

The famous Cauchy determinant identity (see \cite[Thm. 12]{K2}) states that 
$$\det\bigg[\frac{1}{x_i+y_j}\bigg]_{1\le i,j\le n}=\frac{\prod_{1\le i<j\le n}(x_i-x_j)(y_i-y_j)}{\prod_{1\le i,j\le n}(x_i+y_j)}.$$
With the help of this formula, Z.-W. Sun \cite[Thm. 1.4]{S19} proved that for any prime $p\equiv3\pmod4$, 
$$\det\bigg[\frac{1}{i^2+j^2}\bigg]_{1\le i,j\le \frac{p-1}{2}}\equiv \left(\frac{2}{p}\right) \pmod p,$$
where $(\frac{\cdot}{p})$ is the Legendre symbol. In the same paper, Sun also posed the following conjecture.

\begin{conjecture}
Let $p\equiv 2\pmod 3$ be an odd prime and let $$T_p=\left[\f{1}{i^2-ij+j^2}\right]_{1\leq i,j\leq p-1}.$$ Then $2\det T_p$ is a quadratic residue modulo $p$.
\end{conjecture}
Later, the authors and H.-X. Ni \cite{WSN} confirmed this conjecture by proving the $p$-adic congruence 
$$\det T_p\equiv(-1)^{\frac{p+1}{2}}2^{\frac{p-2}{3}}\pmod{p}$$ 
for any odd prime $p\equiv 2\pmod 3$. 

Observe that for any $p$-adic unit $x$, we have the congruence 
$$\frac{1}{x}\equiv x^{p-2}\pmod p.$$
Thus, to obtain a reasonable generalization of $T_p$, it is natural to consider the matrix 
$$D_p(c,d):=[(i^2+cij+dj^2)^{p-2}]_{1\ls i,j\ls p-1},$$
where $c,d$ are integers (note that $T_p=D_p(-1,1)$). 

Z.-W. Sun \cite{S22} initiated the study of $\det D_p(c,d)$. Along this line, recently Luo and Sun \cite{L-S} studied $D_p(1,1)$ where $p\equiv1\pmod3$ and $D_p(2,2)$ where $p\equiv1\pmod4$. For example, they proved that $D_p(2,2)$ is a quadratic residue modulo $p$ whenever $p\equiv 1\pmod8$. 

Motivated by the above results, in this paper, we shall study some arithmetic properties of the determinants of $D_p(c,d)$. 

\subsection{Main Results} We now state our main results. 
We first consider the case $d=1$. It turns out that this case is the most important case (see Theorem \ref{Th1.2} below). For simpilicity, we set
$$D_p(c):=D_p(c,1).$$

\begin{theorem}\label{Th1.1} Let $p$ be an odd prime and let $c$ be an integer. 
Then the following hold. 
\begin{enumerate}
	\item If $c\equiv\pm2\pmod{p}$, then $p\mid\det D_p(c)$.
	\item If $\big(\f{c^2-4}{p}\big)=1$, then there exists $x\in\Z$ such that 
	$$
	\det D_p(c)\equiv x^2\left(\f{-c-2}{p}\right)\pmod{p}.
	$$
	\item If $\big(\f{c^2-4}{p}\big)=-1$, then there exists $x\in\Z$ such that 
	$$
    \det D_p(c)\equiv 2cx^2\left(\f{-c-2}{p}\right)\pmod{p}.
    $$	
\end{enumerate}
\end{theorem}
\begin{remark}
	Part (1) of this theorem is trivial. In fact, $(i^2\pm2ij+j^2)^{p-2}\equiv (i\pm j)^{p-3}\pmod{p}$ for any odd prime $p>3$ and integers $i,j$. An application of \cite[Prop. 9]{K1} will immediately imply the result.
\end{remark}

For general cases, we have the following result.
\begin{theorem} \label{Th1.2}
	Let $p$ be an odd prime.
	Then 
	\begin{enumerate}
		\item $p\mid\det D_p(c,d)$ whenever $\big(\f{d}{p}\big)\not=1$.
		\item 
		$
		\det D_p(c,d^2)\equiv\big(\f{d}{p}\big)\det D_p(cd^{-1})\pmod{p}.
		$
	\end{enumerate}
\end{theorem}

\begin{remark}
	Sun \cite[Thm. 1.3(i)]{S22} proved that $D_p(c,-1)\equiv 0\pmod p$ whenever $p\equiv3\pmod4$. The above result (1) is clearly an extension of Sun's result.  
\end{remark}

As a direct consequence of Theorem \ref*{Th1.1} and Theorem \ref*{Th1.2} (2), we have the following corollary.

\begin{corollary}\label{Corollay. 1}
	Let $p$ be an odd prime and let $c,d$ be integers. Suppose $p\nmid\det D_p(c)$. Then the following results hold.
	
	{\rm (1)} If $\big(\f{c^2-4}{p}\big)=1$, then 
	$$\left(\frac{\det D_p(c)}{p}\right)=\left(\frac{-c-2}{p}\right)^{(p-1)/2}.$$
	
	{\rm (2)} If $\big(\f{c^2-4}{p}\big)=-1$, then 
	$$\left(\frac{\det D_p(c)}{p}\right)=\left(\frac{-c-2}{p}\right)^{(p-1)/2}\left(\frac{2c}{p}\right).$$
	
	{\rm (3)} When $p\nmid d$, we have 
	$$\left(\frac{\det D_p(c,d^2)}{p}\right)=\left(\frac{d}{p}\right)^{(p-1)/2}\left(\frac{\det D_p(cd^{-1})}{p}\right).$$
\end{corollary}

In view of Theorems \ref{Th1.1}--\ref{Th1.2} and Corollary \ref*{Corollay. 1}, what remains to do is to determine when $\det D_{p}(c)$ vanishes modulo $p$, where $c\in\Z$ and $c\not\equiv\pm 2\pmod{p}$. For this, we next introduce the well-known Lucas sequence. 

Given $A,B\in\Z$, the Lucas sequence $u_n(A,B)$ is defined as follows:
\begin{align*}
	&u_0(A,B)=0,\\
	&u_1(A,B)=1,\\
	&u_{n+1}(A,B)=Au_n(A,B)-Bu_{n-1}(A,B).
\end{align*}

By the above notations, we have the following result. 

\begin{theorem}\label{Th1.3}
	Let $p$ be a prime and let $c$ be an integer. Suppose $c\not=\equiv\pm2\pmod{p}$.  Then the following results hold.
	\begin{enumerate}
		\item Suppose $\big(\f{c^2-4}{p}\big)=1$. Then $p\mid\det D_p(c)$ if and only if
		$$
		cu_k(-c,1)+\left[k(4-c^2)+6-c^2\right]u_{k+1}(-c,1)\equiv0\pmod{p}
		$$
		for some $0\ls k\ls p-3$.
		\item Suppose $\big(\f{c^2-4}{p}\big)=-1$. Then $p\mid\det D_p(c)$ if and only if 
		$$
		-cu_k(-c,1)+(2-c^2)u_{k+1}(-c,1)\equiv0\pmod{p}
		$$
		for some $0\ls k\ls p-3$.
	\end{enumerate}
	
\end{theorem}

\subsection{Outline of this paper} The outline of this paper is as follows. We are going to prove some arithmetic properties of trimomial coefficients in Section 2 and our main results will be proved in Section 3.

\section{Some arithmetic properties of trinomial coefficients}
\setcounter{lemma}{0}
\setcounter{theorem}{0}
\setcounter{equation}{0}
\setcounter{conjecture}{0}
\setcounter{remark}{0}
\setcounter{corollary}{0}

For any positive integer $n$ and integer $k$, we define the trinomial coefficients $\binom{n}{k}_c$ by
$$
(x+c+x^{-1})^{n}=\sum_{k=-\infty}^{\infty}\binom{n}{k}_cx^k.
$$
Clearly,  $\binom{n}{k}_c=0$ when $|k|>n$ and for any integer $k$,
\begin{equation} \label{Tri-sym}
	\binom{n}{k}_c=\binom{n}{-k}_c.
\end{equation} By comparing coefficients on both sides of $$(x+c+x^{-1})^n=(x+c+x^{-1})^{n-1}\cdot(x+c+x^{-1}),$$ one can easily verify that 
\begin{equation}\label{Tri-rec}
	\binom{n}{k}_c=\binom{n-1}{k-1}_c+c\binom{n-1}{k}_c+\binom{n-1}{k+1}_c
\end{equation}
In particular, when $k=0$, the following result concerning the generating function of $\binom{n}{0}_c$ can be found in \cite[p. 159]{W}.

\begin{lemma}\label{centri-rec}
	Let $c\in\Z$ with $c^2-4\not=0$. Then
	$$
	\sum_{n=0}^{\infty}\binom{n}{0}_cx^n=\f{1}{\sqrt{1-2cx+(c^2-4)x^2}},
	$$
	which implies that for any positive integer $n$,
	$$
	(n+1)\binom{n+1}{0}_c=(2n+1)c\binom{n}{0}_c-n(c^2-4)\binom{n-1}{0}_c.
	$$
\end{lemma}

\begin{lemma}\label{tri-p} For any integer $k$ and any prime $p$, we have 
	\begin{equation*}
		\binom{p}{k}_c\equiv\begin{cases}
			c\pmod{p}&\t{if $k=0$,}\\
			1\pmod{p}&\t{if $k=\pm p$,}\\
			0\pmod{p}&\t{otherwise.}
		\end{cases}
	\end{equation*}
\end{lemma}
\Proof. Since 
$$
(x+c+x^{-1})^p\equiv x^p+c^p+x^{-p}\equiv x^p+c+x^{-p}\pmod{p},
$$
the conclusion follows from the definition of trinomial coefficients. 
 \qed

For $A,B\in \Z$, the Lucas sequence $u_n(A,B)(n=0,1,2,\ldots)$ is defined as follows:
\begin{align*}
	&u_0(A,B)=0,\\
	&u_1(A,B)=1,\\
	&u_{n+1}(A,B)=Au_n(A,B)-Bu_{n-1}(A,B).
\end{align*}
The following known result can be found in \cite{SZH}.

\begin{lemma}\label{lucas}
	Let $p$ be an odd prime. Suppose $b^2\equiv B\pmod{p}$ where $b\in\Z$. Then
	$$
	u_{\f{p-1}{2}}(A,B)\equiv\begin{cases}
		0\pmod{p}&\t{if $\left(\f{A^2-4B}{p}\right)=1$,}\\ \\
		
		\f{1}{b}\left(\f{A-2b}{p}\right)\pmod{p}&\t{if $\Big(\f{A^2-4B}{p}\Big)=-1$,}
	\end{cases}
    $$
and
	$$
u_{\f{p+1}{2}}(A,B)\equiv\begin{cases}
	\left(\f{A-2b}{p}\right)\pmod{p}&\t{if $\left(\f{A^2-4B}{p}\right)=1$.}\\ \\ 
	0\pmod{p}&\t{if $\left(\f{A^2-4B}{p}\right)=-1$,}
    \end{cases}
    $$
\end{lemma}
We also need the following lemma which appeared in \cite[(2.5)]{S14}.

\begin{lemma} \label{cen-tri}
	Let $c$ be an integer. For any odd prime $p$, we have
	$$
	\binom{p-1}{0}_c\equiv\left(\f{c^2-4}{p}\right)\pmod{p}.
	$$
\end{lemma}

By Lemma \ref{lucas} and Lemma \ref{cen-tri}, we are now ready to prove the following result.

\begin{lemma}\label{tri-p-2}
	Let $p$ be an odd prime $p$. For any integer $c$ with $c\not\equiv\pm2\pmod{p}$, we have
$$
\binom{p-2}{0}_c\equiv\f{c}{c^2-4}\left(\f{c^2-4}{p}\right)\pmod{p}
$$
and
\begin{equation*}
\binom{p-2}{\f{p-1}{2}}_c\equiv\begin{cases}
	\f{c}{2(c^2-4)}\big(\f{-c-2}{p}\big)\pmod{p}&\t{if $\big(\f{c^2-4}{p}\big)=1$,}\\
	-\f{1}{c^2-4}\big(\f{-c-2}{p}\big)\pmod{p}&\t{if $\big(\f{c^2-4}{p}\big)=-1$.}
\end{cases}	
\end{equation*}

\end{lemma}

\Proof. By Lemma \ref{centri-rec}, we have
$$
p\binom{p}{0}_c=(2p-1)c\binom{p-1}{0}_c-(p-1)(c^2-4)\binom{p-2}{0}_c,
$$
and hence by Lemma \ref{cen-tri}, we obtain that
$$
\binom{p-2}{0}_c\equiv\f{c}{c^2-4}\left(\f{c^2-4}{p}\right)\pmod{p}
$$
as desired. 

Combining \eqref{Tri-rec} with Lemma \ref{tri-p}, we obtain that
\begin{equation}\label{p-1-rec}
	\binom{p-1}{k+1}_c\equiv-c\binom{p-1}{k}_c-\binom{p-1}{k-1}_c\pmod{p},\ 1\ls k\ls p-1,
\end{equation}
and 
\begin{equation}\label{p-1-(-1,0,1)}
	\binom{p-1}{1}_c+c\binom{p-1}{0}_c+\binom{p-1}{-1}_c\equiv c\pmod{p}.
\end{equation}

It follows from \eqref{Tri-sym}, \eqref{p-1-(-1,0,1)} and Lemma \ref{cen-tri} that
\begin{equation}\label{p-1-1}
	\binom{p-1}{1}_c\equiv\f{1-\big(\f{c^2-4}{p}\big)}{2}c\pmod{p},
\end{equation}
and hence by \eqref{p-1-rec}, when $0\ls k\ls p$, we have
\begin{equation}\label{tri-p-1}
\binom{p-1}{k}_c\equiv\f{1+\big(\f{c^2-4}{p}\big)}{2}cu_k(-c,1)+\Big(\f{c^2-4}{p}\Big)u_{k+1}(-c,1)\pmod{p}.
\end{equation}
Combining this with Lemma \ref{lucas}, we get
\begin{align}\label{p-1-(p-3)/2}
\binom{p-1}{\f{p-3}{2}}_c\equiv
\begin{cases}
-c\big(\f{-c-2}{p}\big)\pmod{p}&\mbox{if $\big(\f{c^2-4}{p}\big)=1$,}\\
-\big(\f{-c-2}{p}\big)\pmod{p} &\mbox{if $\big(\f{c^2-4}{p}\big)=-1$,}
\end{cases}
\end{align}
and
\begin{align}\label{p-1-(p+1)/2}
\binom{p-1}{\f{p+1}{2}}_c\equiv
\begin{cases}
0\pmod{p}                     &\mbox{if $\big(\f{c^2-4}{p}\big)=1$,}\\
\big(\f{-c-2}{p}\big)\pmod{p} &\mbox{if $\big(\f{c^2-4}{p}\big)=-1$.}
\end{cases}
\end{align}

By definitions of trinomial coeffiticents, we have
$$
(x+c+x^{-1})^{p-1}=\sum_{k=-p+1}^{p-1}\binom{p-1}{k}_cx^k.
$$
Taking derivative on both sides and we get
$$
(p-1)(x+c+x^{-1})^{p-2}(1-\f{1}{x^2})=\sum_{k=-p+1}^{p-1}k\binom{p-1}{k}_cx^{k-1}.
$$
Comparing the coefficient of $x^{k-1}$ on both sides and we obtain that
\begin{equation*}
k\binom{p-1}{k}_c=(p-1)\binom{p-2}{k-1}_c-(p-1)\binom{p-2}{k+1}_c.
\end{equation*}
Combining this with \eqref{Tri-rec}, we have
\begin{equation}\label{p-1-p-2}
\begin{aligned}
&(k+1)\binom{p-1}{k-1}_c-(k-1)\binom{p-1}{k+1}_c\\
\equiv&2\binom{p-1}{k-1}_c+2\binom{p-1}{k+1}_c+\left(\binom{p-2}{k}_c-\binom{p-2}{k-2}_c\right)-\left(\binom{p-2}{k+1}_c-\binom{p-2}{k}_c\right)\\
\equiv&2\binom{p-1}{k-1}_c+2\binom{p-1}{k+1}_c+2\binom{p-2}{k}_c-\left(\binom{p-1}{k-1}_c-c\binom{p-2}{k-1}_c-\binom{p-2}{k}_c\right)\\
      &-\left(\binom{p-1}{k+1}_c-c\binom{p-2}{k+1}_c-\binom{p-2}{k}_c\right)\\
\equiv&\binom{p-1}{k-1}_c+\binom{p-1}{k+1}_c+c\binom{p-2}{k-1}_c+4\binom{p-2}{k}_c+c\binom{p-2}{k+1}_c\\
\equiv&\binom{p-1}{k-1}_c+c\binom{p-1}{k}_c+\binom{p-1}{k+1}_c+(4-c^2)\binom{p-2}{k}_c\\
\equiv&\binom{p}{k}_c+(4-c^2)\binom{p-2}{k}_c\pmod{p}.\\
\end{aligned}
\end{equation}
With the aid of this, Lemma \ref{tri-p}, \eqref{p-1-(p-3)/2} and \eqref{p-1-(p+1)/2}, via simple computation,  the reader may derive the desired valuation of $\binom{p-2}{\f{p-1}{2}}_c$.
\qed

\section{Proofs of The Main Results}
We first introduce the definition of circulant matrices. Let $R$ be a commutative ring. Let $m$ be a positive integer and $a_0,a_1,\ldots,a_{m-1}\in R$. We define the circulant matrix $C(a_0,\ldots,a_{m-1})$ to be an $m\times m$ matrix whose ($i$-$j$)-entry is $a_{i-j}$ where the indices are cyclic module $m$. 
The second author \cite[Lemma 3.4]{W21}  proved the following lemma.

\begin{lemma} \label{lemma-ffa}
	Let $R$ be a commutative ring. Let $m$ be a positive integer. Let $a_0,a_1,\ldots,a_{m-1}\in R$ such that
	\begin{align*}
	a_i=a_{m-i} \ \ \t{for each $1\ls i\ls m-1$.}
	\end{align*}
If $m$ is even, then there exists an element $u\in R$ such that 
$$
\det C(a_0,\ldots,a_{m-1})=\left(\sum_{i=0}^{m-1}a_i\right)\left(\sum_{i=0}^{m-1}(-1)^ia_i\right)\cdot u^2.
$$
If $m$ is odd, then there exists an element $v\in R$ such that
$$
\det C(a_0,\ldots,a_{m-1})=\left(\sum_{i=0}^{m-1}a_i\right)\cdot v^2.
$$
\end{lemma}

\medskip
\noindent{\bf Proof of Theorem \ref{Th1.1}}.
Suppose that $g$ is a primitive element modulo $p$. Then 
\begin{align*}
\det D_p(c)\equiv&\det\left[(g^{2i}+cg^{i+j}+g^{2j})^{p-2}\right]_{0\ls i,j\ls p-2}\\ \equiv&\det\left[g^{i-j}(g^{2(i-j)}+cg^{i-j}+1)^{p-2}\right]_{0\ls i,j\ls p-2}\pmod{p}.
\end{align*}
For $0\ls i\ls p-2$, set
$$
a_{i}=g^i(g^{2i}+cg^{i}+1)^{p-2}.
$$
Note that
\begin{align*}
	a_{p-1-i}=g^{p-1-i}(g^{2(p-1-i)}+cg^{p-1-i}+1)^{p-2}\equiv g^{i}(1+cg^i+g^{2i})^{p-2}\equiv a_i\pmod{p}.
\end{align*}
 Hence by Lemma \ref{lemma-ffa}, there exists $u\in\Z$ such that 
\begin{equation}
	\det D_p(c)\equiv\left(\sum_{i=0}^{p-2}a_i\right)\left(\sum_{i=0}^{p-2}(-1)^ia_i\right)\cdot u^2\pmod{p}.
\end{equation}
Furthermore, we have
\begin{align*}
	\sum_{i=0}^{p-2}a_i&\equiv\sum_{i=1}^{p-1}i(i^2+ci+1)^{p-2}\\
	                   &\equiv\sum_{i=1}^{p-1}\sum_{k=-p+2}^{p-2}\binom{p-2}{k}_ci^k\\
	                   &\equiv\sum_{k=-p+2}^{p-2}\sum_{i=1}^{p-1}\binom{p-2}{k}_ci^k\\
	                   &\equiv-\binom{p-2}{0}_c\pmod{p}
\end{align*}
and
\begin{align*} 
	\sum_{i=0}^{p-2}(-1)^ia_i&\equiv\sum_{i=1}^{p-1}\Big(\f{i}{p}\Big)i(i^2+ci+1)^{p-2}\\
	&\equiv\sum_{i=1}^{p-1}i^{\f{p-1}{2}}(i+c+i^{-1})^{p-2}\\
	&\equiv\sum_{k=-p+2}^{p-2}\sum_{i=1}^{p-1}\binom{p-2}{k}_ci^{k+\f{p-1}{2}}\\
	&\equiv-\binom{p-2}{-\f{p-1}{2}}_c-\binom{p-2}{\f{p-1}{2}}_c\\
	&\equiv-2\binom{p-2}{\f{p-1}{2}}_c\pmod{p}.
\end{align*}
Therefore,
$$
\det D_p(c)=2\binom{p-2}{0}_c\binom{p-2}{\f{p-1}{2}}_c\cdot u^2,
$$
and the conclusion follows from Lemma \ref{tri-p-2}. This completes the proof.\qed

\medskip

We next turn to the proof of Theorem \ref{Th1.2}. 
Fix any nonzero element $a\in\mathbb{F}_p$, the map $x\mapsto ax$ induces a permutatin $\pi_a$ of $\mathbb{F}_p$. Lerch \cite{L} computed the sign of this permutation explicitly.

\begin{lemma}\label{sign}
	Assume that $p$ is an odd prime. Let $\sgn(\pi_a)$ denote the sign of the permutation $\pi_a$. Then
	\begin{equation*}
		\sgn(\pi_a)=\Big(\f{a}{p}\Big).
	\end{equation*}
\end{lemma}
Now we are in a position to prove our second result. 

\noindent{\bf Proof of Theorem \ref{Th1.2}}.
By Lemma \ref{sign}, we have 
\begin{align*}
	\det D_p(c,d)&=\det\left[(i^2+cij+dj^2)^{p-2}\right]_{1\ls i,j\ls p-1}\\
	&\equiv\det\left[(di^2+cdij+(dj)^2)^{p-2}\right]_{1\ls i,j\ls p-1}\\
	&\equiv\Big(\f{d}{p}\Big)\det\left[(di^2+cij+j^2)^{p-2}\right]_{1\ls i,j\ls p-1}\\
	&\equiv\Big(\f{d}{p}\Big)\det D_p(c,d).
\end{align*}
Hence $\det\left[(i^2+cij+dj^2)^{p-2}\right]_{1\ls i,j\ls p-1}\equiv0\pmod{p}$  whenever $\big(\f{d}{p}\big)\neq1$.

In addition, by Lemma \ref{sign} again, 
\begin{align*}
	\det D_p(c,d^2)&=\det\left[(i^2+cd^{-1}i(dj)+(dj)^2)^{p-2}\right]_{1\ls i,j\ls p-1}\\
	&=\Big(\f{d}{p}\Big)\det\left[(i^2+cd^{-1}ij+j^2)^{p-2}\right]_{1\ls i,j\ls p-1}.\\
	&\equiv\Big(\f{d}{p}\Big)D_p\left(cd^{-1}\right)\pmod{p}.
\end{align*}

This completes our proof. \qed

To determine whether $\det D_p(c)\mod p$ vanishes, we need the following result (see  \cite[Lemma 10]{K2}).
\begin{lemma}\label{p(xy)}
	Let $P(Z)=a_{n-1}Z^{n-1}+a_{n-2}Z^{n-2}+\cdots+a_0$ be a poynomial over $\mathbb{F}_p$. Then
	$$
	\det[P(X_iY_j)]_{1\ls i.j\ls n}=\prod_{i=0}^{n-1}a_i\prod_{1\ls i<j\ls n}(X_i-X_j)(Y_i-Y_j).
	$$
\end{lemma}
Now we are in a position to prove our last result.

\noindent{\bf Proof of Theorem \ref{Th1.3}}. Note that
\begin{align*}
\det D_p(c)&=\prod_{j=1}^{p-1}j^{2p-4}\cdot\det\left[\left(\left(\f{i}{j}\right)^2+c\f{i}{j}+1\right)^{p-2}\right]_{1\ls i,j\ls p-1}\\
&\equiv\det\left[\sum_{k=-p+2}^{p-2}\binom{p-2}{k}_c\left(\f{i}{j}\right)^{k+p-2}\right]_{1\ls i,j\ls p-1}\\
&\equiv\det\left[\sum_{k=0}^{p-2}\left(\binom{p-2}{p-2-k}_c+\binom{p-2}{k+1}_c\right)\left(\f{i}{j}\right)^{k+p-2}\right]_{1\ls i,j\ls p-1}\pmod{p}.\\
\end{align*}
Hence an application of Lemma \ref{p(xy)} shows that $\det D_p(c)\equiv0\pmod{p}$ if and only if 
$$
\binom{p-2}{p-2-k}_c+\binom{p-2}{k+1}_c\equiv0\pmod{p}
$$
for some $0\ls k\ls p-2.$ 

Note that when $k=p-2$, $\binom{p-2}{0}_c+\binom{p-2}{p-1}_c=\big(\f{c^2-4}{p}\big)\not=0$.
If $0\ls k\ls p-3$, by \eqref{p-1-p-2}, we have
\begin{align*}
&(4-c^2)\left[\binom{p-2}{p-2-k}_c+\binom{p-2}{k+1}_c\right]\\
\equiv &(k+3)\binom{p-1}{p-1-k}_c+(k+2)\binom{p-1}{k}_c-(k+1)\binom{p-1}{p-3-k}_c-k\binom{p-1}{k+2}_c\pmod{p}.
\end{align*}
By \cite[Lemma 4.2]{L-S}, $\binom{p-1}{p-k}_c\equiv u_k(-c,1)\pmod{p}$. Combining this with \eqref{tri-p-1}, the reader may verify that
\begin{align*}
&(4-c^2)\left[\binom{p-2}{p-2-k}_c+\binom{p-2}{k+1}_c\right]\\
\equiv &c\Big(\f{c^2-4}{p}\Big)u_k(-c,1)+\left[2k+4-\f{k+2}{2}c^2+\left(2k+2-\f{k}{2}c^2\right)\Big(\f{c^2-4}{p}\Big)\right]u_{k+1}(-c,1)\pmod{p}.
\end{align*}

In view of the above, we have completed the proof. \qed

\Ack\ This work was supported by the Natural Science Foundation of China (Grant no. 12101321).

\end{document}